\documentclass[11pt]{amsart}

\usepackage{times,amsmath,amsthm,amsfonts,amssymb}
\usepackage{mathrsfs}
\usepackage{fullpage}

\title[Interpolation and Sampling Hypersurfaces on the Ball]{Sufficient Conditions for Interpolation and Sampling Hypersurfaces in the Bergman Ball}

\author{Tam\' as Forg\' acs}
\address{Department of Mathematics\newline \indent
University of Illinois, Urbana, IL 61801}
\email{forgacs@uiuc.edu}

\author{Dror Varolin}
\address{Department of Mathematics\newline\indent
Stony Brook University, Stony Brook, NY 11794}
\email{{\tt dror@math.sunysb.edu}}

\thanks{Research partially supported by NSF grant DMV-0400909}
\date{}

\newcommand{\noi}{\noindent}

\newcommand{\cc}{{\mathcal C}}

\newcommand{\co}{{\mathcal O}}

\newcommand{\sC}{{\mathscr C}}
\newcommand{\sd}{{\mathscr D}}

\newcommand{\sh}{{\mathscr H}}

\newcommand{\sL}{{\mathscr L}}

\newcommand{\so}{{\mathscr O}}
\newcommand{\sP}{{\mathscr P}}

\newcommand{\sw}{{\mathscr W}}

\newcommand{\fH}{{\mathfrak H}}

\newcommand{\vp}{\varphi}
\newcommand{\ve}{\varepsilon}

\newcommand{\C}{{\mathbb C}}
\newcommand{\D}{{\mathbb D}}

\newcommand{\R}{{\mathbb R}}

\newcommand{\di}{\partial}
\newcommand{\dbar}{\bar \partial}
\newcommand{\ii}{\sqrt{-1}}

\newcommand{\re}{{\rm Re\ }}

\newcommand{\relcomp}{\subset \subset}

\begin{document}

\maketitle

\theoremstyle{plain}
\newtheorem{thm}{\sc Theorem}[section]
\newtheorem{lem}[thm]{\sc Lemma}
\newtheorem{o-thm}[thm]{\sc Theorem}
\newtheorem{prop}[thm]{\sc Proposition}
\newtheorem{cor}[thm]{\sc Corollary}

\theoremstyle{definition}
\newtheorem{conj}[thm]{\sc Conjecture}
\newtheorem{defn}[thm]{\sc Definition}
\newtheorem{qn}[thm]{\sc Question}
\newtheorem{ex}[thm]{\sc Example}
\newtheorem{pr}[thm]{\sc Problem}

\theoremstyle{remark}
\newtheorem*{rmk}{\sc Remark}
\newtheorem*{ack}{\sc Acknowledgment}

\section{Introduction}

Recall that the Bergman metric on the unit ball $B = \left \{ z \in \C ^n \ ;\ |z| < 1 \right \}$ is the K\" ahler metric whose associated $(1,1)$-form is $\omega _B = - (n+1)d d^c \lambda,$ where
\[
\lambda = \log (1-|z|^2) - \frac{n}{n+1} \log (n+1)
\]
and in our convention $d^c= \frac{\ii}{2}(\dbar -\di).$  The weighted Bergman spaces on the Bergman ball are
\[
\sh^2 (B,\kappa) := \left \{ F \in \so (B)\ ;\ \int _B |F|^2
e^{-\kappa} \omega _B ^n < +\infty \right \},
\]
where $\so (X)$ denotes the space of holomorphic functions on a complex manifold $X$.  In this paper we assume that $\kappa$ is $\sC ^2$.  The case $\kappa = -(n+1) \log (1-|z|^2)$ corresponds to the classical Bergman space of holomorphic functions that are square integrable with respect to Lebesgue measure.

Given a nonsingular closed complex hypersurface $W \subset B$, we let
\[
\fH ^2 (W,\kappa) := \left \{ f \in \so (W)\ ;\ \int _W |f|^2 e^{-\kappa} \omega _B ^{n-1} < +\infty \right \}.
\]
\begin{defn}
\begin{enumerate}
\item[(a)] We say that $W$ is an interpolation hypersurface if for each $f \in \fH ^2 (W,\kappa)$ there exists $F \in \sh ^2(B,\kappa)$ such that $F|W=f$.

\item[(b)] We say that $W$ is a sampling hypersurface if there is a constant $A$ such that for every $F \in \sh ^2 (B,\kappa)$,
\begin{eqnarray}\label{samp-ineq}
\frac{1}{A} \int _B |F|^2 e^{-\kappa} \omega ^n _B \le \int _W |F|^2 e^{-\kappa} \omega ^{n-1} _B \le A \int _B |F|^2 e^{-\kappa} \omega ^n _B.
\end{eqnarray}
\end{enumerate}
\end{defn}

\noi  Let $F_a$ denote a holomorphic involution of $B$ sending $0$ to $a$ (see Section \ref{bg-rev}).

\begin{rmk}
We will often use, without explicit indication, the fact that $F_z$ is an involution.  Thus the reader should not be confused if $F_z$ is seen when $F_z^{-1}$ is expected.
\end{rmk}

We define the {\it total density tensor} of $W$ in the ball of radius $r$ to be the $(1,1)$-form
\[
\Upsilon _r ^W (z) = \frac{1}{V_n(r)}\left (  \int_{B(0,r)} \frac{\di ^2 \log |T(F_{z}(\zeta))|^2}{\di z ^i \di \bar z^j}
\omega _B ^n \right )  \ii  dz^i \wedge d\bar z ^j.
\]
Here $T$ is any holomorphic function such that $W = \{ T = 0 \}$
with $dT|W$ nowhere zero, and
\[
V_n(r) = \int _{B(0,r)} \omega _B ^n
\]
is the volume of the Euclidean ball of radius $r$ and center $0$, with respect to (our normalization of) the volume induced by the Bergman metric.  The total density tensor is a Bergman ball analog of the total density tensor introduced in \cite{osv} in the case of $\C^n$.  In the case of the Bergman ball, some of the more basic properties of the total density tensor do not follow as readily as their analogs in the $\C ^n$ case.  For example, at the end of Section \ref{bg-rev} we will show that the definition of $\Upsilon ^W _r$ is independent of the choice of $T$.

We define
\[
[W]_{\ve}(z)  = dd^c \left ( \frac{1}{V_n(\ve)} \int _{|F_z(\zeta)|< \ve} \log |T(\zeta)|^2 \omega _B^n (\zeta) \right ).
\]
If we denote by $[W]$ the current of integration along $W$, then
$[W]_{\ve}$ is in some sense the average of $[W]$ over the
Bergman-Green ball of radius $\ve$.  Note that $[W]_{\ve}=
\Upsilon ^W_{\ve}$ and thus $[W]_{\ve}$ is independent of the choice
of $T$.   Moreover, though not necessarily smooth, the current
$[W]_{\ve}$ is locally bounded, as can be seen by changing
variables in the intergral and then differentiating under the
integral. Finally, it is also clear that, in the sense of
currents, $[W]_{\ve} \to [W]$ as $\ve \to 0$.

\begin{defn}\label{densities}\
\begin{enumerate}
\item[(I)]  Let $\sP _W (B)$ denote the set of $(n-1,n-1)$-forms $\theta$ on $B$ with the following properties.
\begin{enumerate}
\item $\theta \wedge \omega _B \ge c e^{n\lambda} \omega _B ^n$ for some constant $c > 0$.

\item For each $\ve > 0$ there exists $C>0$ such that $[W] _{\ve} \wedge \theta \le C [W]_{\ve} \wedge \omega _B ^{n-1}$.

\item $dd^c \theta = 0$.
\end{enumerate}

\item[(II)]  For $\theta \in \sP _W(B)$, let
$$
\sd _B^+(W,\kappa)[\theta] = \limsup _{r \to 1} \sup _{z\in B} \frac{\left ( \Upsilon _r ^W  + \frac{n}{n+1} \omega _B \right ) \wedge \theta }{\ii \di \dbar \kappa \wedge \theta}
$$
and
$$
\sd _B^-(W,\kappa)[\theta] = \liminf _{r \to 1} \inf _{z\in B} \frac{\left ( \Upsilon _r ^W  + \frac{n}{n+1} \omega _B \right ) \wedge \theta }{\ii \di \dbar \kappa \wedge \theta}
$$

\item[(III)]  The upper and lower densities of $W$ are
$$
\sd _B^+(W,\kappa) = \sup _{\theta \in \sP _W(B)} \sd _B^+(W,\kappa)[\theta]
$$
and
$$
\sd _B^-(W,\kappa) = \sup _{\theta \in \sP _W(B)} \sd_B^- (W, \kappa ) [\theta]
$$
\end{enumerate}
\end{defn}

\noi From here on out we assume that $W$ is uniformly flat (see Section \ref{uf-section} for the definition) and that
\[
\frac{1}{C} \omega _B \le \ii \di \dbar \kappa \le C \omega _B
\]
for some constant $C > 1$.  Our main results are the following two theorems.

\begin{thm}\label{interp-thm}
If $\sd ^+_B (W,\kappa) < 1$, then $W$ is an interpolation hypersurface.
\end{thm}

\begin{thm}\label{samp-thm}
If $\sd ^-_B (W,\kappa ) > 1$, then $W$ is a sampling hypersurface.
\end{thm}

Theorems \ref{interp-thm} and \ref{samp-thm} give generalizations to higher dimensions of results of Seip \cite{sp} and of Berndtsson-Ortega Cerd\`a \cite{quimbo}.  By now Theorems \ref{interp-thm} and \ref{samp-thm} carry with them a rich history.  Most recently, results analogous to Theorems \ref{interp-thm} and \ref{samp-thm} have been established for the case of $\C ^n$ in the paper \cite{osv}, which we refer to for further historical remarks regarding interpolation and sampling problems for Bergman spaces.

Though there is a strong similarity between the results of \cite{osv} and the present paper, the methods of proof are completely different.  In fact, the present approach and the approach of \cite{osv} could be used interchangeably for the case of $\C ^n$ and the Bergman ball.

In the case of interpolation, we employ the Ohsawa-Takegoshi technique to extend functions from the submanifold $W$ to the ball in one shot, rather than using the $L^2$ Cousin I-type approach to extend the function locally and then patch together the resulting local extensions.  (We should perhaps remark that if one wants to apply the Cousin I-type method in the case of the ball, then the negativity of the curvature of $\omega _B$ requires the use of a sharper version of H\" ormander's $\dbar$ Theorem, due to Ohsawa.  The need for Ohsawa's Theorem was already noticed in the 1-dimensional case \cite{quimbo}.)

By contrast with \cite{osv}, our approach to sampling is closer in spirit to the technique that has been used in the one-variable case in \cite{quimbo}.  Our densities, laid out in Definition \ref{densities} above, do not directly correspond to those in \cite{osv} (though we prove in Section \ref{dens-again} that they are actually the same).  We feel that the methods of the present paper fit in more naturally with the Hilbert Space approach.  The proofs also seem more elementary than the Beurling-inspired approach used in \cite{osv}.

The paper is organized as follows.

\tableofcontents


\begin{ack}
Thanks to Jeff McNeal and Quim Ortega-Cerd\` a for many stimulating discussions.
\end{ack}

\section{Rapid review of Bergman geometry}\label{bg-rev}

Bergman geometry is one of the oldest and most studied areas of complex geometry.  Therefore we content ourselves with stating facts, and provide few proofs.

\subsubsection*{Bergman metric}

As already mentioned, the Bergman metric is $\omega _B = -(n+1)dd^c \lambda.$ It is easy to see that, with $\omega _E = dd^c |z|^2$ denoting the Euclidean metric,
$$
\left . \omega _B\right  |_{z=0} = (n+1) \left . \omega _E \right
|_{z=0} \quad {\rm and} \quad \omega _B ^n = e^{-(n+1)\lambda } \omega _E ^n,
$$
and in particular,
$$
{\rm Ricci} (\omega _B) = - \omega _B.
$$

\subsubsection*{Basics of Aut(B)}

For the reader's convenience, we recall that Aut$(B)$ contains the involutions
$$
F_a(z) = \frac{a-P_az-s_aQ_az}{1-\langle z,a\rangle}, \quad a\in
B - \{ 0 \}, \quad F_0(z) = -z,
$$
where $P_a = |a|^{-2} aa^{\dagger}$, $Q_a = I-P_a$ and $s_a =
\sqrt{1-|a|^2}$. Moreover, the Schwarz Lemma shows that any automorphism of $B$ is of the form $ UF_a$ for some unitary $U$.  Note that $F_a(0)=a$  and
$$
1-|F_a(z)|^2 = \frac{(1-|z|^2)(1-|a|^2)}{|1-\left < z , a \right >|^2}.
$$
Thus ${\rm Aut}(B)$ acts transitively on the ball and $\omega _B$ is Aut$(B)$-invariant.  (For much more detail on this and the next paragraph, the reader is referred to \cite{rudin} or \cite{st}.)

\subsubsection*{Basic potential theory of the Bergman metric}

Recall that the Bergman Laplacian $\Delta _B$ associated to
$\omega _B$ is the $\omega _B$-trace of $dd^c$:
$$
\left (\Delta _B g \right ) \omega _B ^n= dd^c g \wedge \omega _B ^{n-1}.
$$
\begin{defn}
The Green's function with pole at $a\in B$ is the function $G_B (z,a)$ satisfying
$$
\Delta _B (G_B(\cdot,a)) \omega _B ^n = \delta _a \quad and \quad
G_B(\cdot,a)|\di B = 0.
$$
\end{defn}

\noi Using Aut$(B)$-invariance, it is easily seen that $G(z,a) =
G(F_a(z), 0)$ and that
$$
n(n+1)\left ( \Delta _B g \right )(a) = {\rm trace} \left ( DF_a(0) ^{\dagger} D^{1,1}g(a) DF_a(0) \right ) .
$$
Here $D^{1,1}g$ is the matrix of the $(1,1)$-form $\ii \di \dbar
g$ in Euclidean coordinates. Setting $\gamma_B = G_B (\cdot , 0),$ we see from unitary invariance that $\gamma _B (z) = f(|z|^2)$ for some function $f$. Substitution into the Bergman-Laplace equation and solving the resulting ODE shows that
$$
f(t) = - C_n \int _t ^1\frac{(1-u)^{n-1}}{u^n} du,
$$
where $C_n =(2\pi)^{-n}(n+1)^{-(n-1)}.$

Note that $f'(t) > 0$.  It follows that for each $a \in B$ the sublevel sets $G(z,a)$ are also the sublevel sets of $|F_a(z)|$.  We use the latter to define distances.

\begin{defn}
{\rm (i)}  The Bergman-Green distance between two points $a$ and $b$ in $B$ is
\[
|F_a(b)|.
\]
{\rm (ii)}  The Bergman-Green balls with center $a$ and radius $r$ are
\[
E(a,r) = F_a(B(0,r)) = \{ z\in B\ ;\ |F_a(z)|<r \}.
\]
\end{defn}

By using the Green-Stokes identity
\begin{eqnarray}\label{sg-id}
\int _{\di D}(g_1 d^c g_2 - g_2 d^c g_1) \wedge \omega ^{n-1} =\int _{D}(g_1 dd^c g_2 - g_2 dd^c g_1)\wedge \omega ^{n-1},
\end{eqnarray}
where $(D,\omega)$ is an $n$ dimensional K\" ahler manifold with boundary and $g_1, g_2 : D \to \C$ are functions, we obtain the following Lemma.
\begin{lem}\label{mvp-lem}
Let $h$ be a function such that $\Delta _B h \ge 0$.  Then
\begin{eqnarray}\label{mvp}
h(0) \le  \frac{1}{(2\pi)^n} \int _{\di B} h(rz) d^c |z|^2 \wedge
\omega _E ^{n-1}(z).
\end{eqnarray}
Moreover, equality holds when $\Delta _B h \equiv 0$.
\end{lem}

\begin{proof}
Apply (\ref{sg-id}) with $D=B(0,r)$, $r < 1$, $g_1 = h$ and $g_2 = \gamma _r$, where
$$
\gamma _r (z) := \gamma _B (z) + C_n \int _{r^2} ^1
\frac{(1-t)^{n-1}}{t^n} dt,
$$
observing that $\gamma _r |\di B(0,r) \equiv 0$ and $dd^c \gamma _r \wedge \omega _B ^{n-1} = \delta _0.$  The result now follows by direct computation.
\end{proof}

\begin{cor}\label{berg-ineq}
Let $h$ be a function such that $\Delta _Bh\ge 0$.  Then for all $r<1$,
\begin{eqnarray}\label{mvp-berg}
h(0) \le  \frac{1}{V_n(r)} \int _{B(0,r)} h \omega _B ^n.
\end{eqnarray}
Moreover, equality holds in {\rm (\ref{mvp-berg})} when $\Delta _B h \equiv 0$.
\end{cor}

Let us end this section by justifying our claim that $\Upsilon ^W _r$ is independent of the choice of holomorphic function $T$ defining $W$.  Suppose $\tilde T$ is another function such that $W = \{ \tilde T=0\}$ and $d\tilde T |W$ is free of zeros.  Then the function $\tilde T / T$ is holomorphic and free of zeros in the ball.  Since the ball is simply connected, any zero-free holomorphic function is the exponential of some holomorphic function.  Thus $\tilde T = e^h T$ for some holomorphic function $h$.  It follows that
\[
\int _{B(0,r)} \log |\tilde T (F_z (\zeta))|^2 \omega _B ^n = \int _{B(0,r)} \log | T (F_z (\zeta))|^2 \omega _B ^n + 2  \int _{B(0,r)} \re h (F_z (\zeta)) \omega _B ^n.
\]
Since $\re h$ is (pluri)harmonic, its ball average, with respect to a radially symmetric probability measure, is equal to its central value.  Since $F_z (0)=z$, we have
\[
\int _{B(0,r)} \log |\tilde T (F_z (\zeta))|^2 \omega _B ^n = \int _{B(0,r)} \log | T (F_z (\zeta))|^2 \omega _B ^n + 2 V_n(r) \re h (z).
\]
The pluriharmonicity of $\re h$ thus completes the justification of our claim.

\section{Uniform flatness}\label{uf-section}

In \cite{osv} a notion of uniform flatness was developed for
closed smooth hypersurfaces in $\C ^n$.  Here we define the
analogous notion for the ball with its Bergman geometry.

Let
$$
N_{\ve}^B(W)=\left \{ z \in B \ \big| \ \inf_{w \in W} |F_z(w)| < \ve \right \}
$$
\begin{defn}\label{uni-flat}
We define a smooth divisor $W$ in $B$ to be uniformly flat if there exists an $\ve _0 >0$ such that $N_{\ve _0}^B(W)$ has the following property:  for each $z \in N_{\ve _0}(W)$ there is a unique $w_z \in W$ minimizing the "distance to $z$" function $w \mapsto |F_z(w)|$ along $W$.
\end{defn}

\begin{rmk}
Recall that a pseudohyperbolic disk of radius $\ve$ is the image under $F \in {\rm Aut} (B)$ of the disk $\{ (0,...,0,z) \in B\ ;\ |z| < \ve\}$.  The unifrom flatness hypothesis implies that in fact $N_{\ve _o}^B (W)$ is foliated by pseudohyperbolic disks.  Indeed, since our condition is invariant under ${\rm Aut}(B)$, it suffices to see this for the case where $W \ni 0$ and $T_{W,0} = \{ z_n = 0\}$. In this case, it is clear that the boundary of the disk $\{ (0,...,0,z) \in B\ ;\ |z| < \ve _o\}$ has distance exactly $\ve _o$ to the origin.

These observations imply the existence of a diffeomorphism
\[
\Phi  : W \times \D (0,\ve_o) \to N^B_{\ve _o}(W)
\]
that is holomorphic in the disk variable, and sends each disk $\{w\} \times \D (0,\ve _o)$ to the disk with center at $w$, which minimizes the pseudo-hyperbolic distance and whose tangent vector is orthogonal to $T_{W,w}$ in the Bergman metric.
\end{rmk}

\noi The following consequence of uniform flatness is useful.

\begin{lem}\label{u-f-char}
If a closed non-singular complex hypersurface $W \subset B$ is uniformly flat, then
there exist $\ve _0 > 0$ and $C>0$ such that for each $z \in W$ the set $F_z (W) \cap B(0,\ve _0)$ is a graph, over the Euclidean $\ve _0$-neighborhood of the origin in the tangent space $T_{F_z(W), 0} = dF_z(T_{W,z})$, of some function $f$ such that
$$
|f(x)| \le C |x|^2, \quad |x| < \ve _0.
$$
\end{lem}

\begin{proof}[Sketch of proof]
Since the notion of uniform flatness is invariant with respect to ${\rm Aut}(B)$, it suffices to assume that $z = 0 \in W$.  Moreover, since we are working in a small neighborhood, we may replace the Bergman metric by the Euclidean metric, and the pseudo-hyperbolic distance $|F_z(w)|$ by Euclidean distance.    In this setting, the result follows from Proposition 2.2 in \cite{osv}.  We leave the details to the interested reader.
\end{proof}

\section{The density conditions again}\label{dens-again}

\subsection{Reformulation of the density conditions}

It will be useful to rewrite the hypotheses on the upper and lower densities in terms of the positivity properties of certain associated differential forms.
\begin{lem}\label{density-conditions}
\begin{enumerate}
\item[1.] If $\sd ^+_B (W,\kappa) < 1$, then there is a positive constant $c$ such that
\[
\ii \di \dbar \kappa - \frac{n}{n+1} \omega _B - \Upsilon ^W _r
\ge c \ii \di \dbar \kappa.
\]

\item[2.] If $\sd ^-_B(W,\kappa) > 1$, then there exists $\theta \in \sP _W(B)$ and $c > 0$ such that
$$
\left ( \Upsilon ^W _r + \frac{n}{n+1} \omega _B - \ii \di \dbar
\kappa \right ) \wedge \theta \ge c e^{n\lambda} \omega _B ^n.
$$
\end{enumerate}
\end{lem}

\begin{proof}
After using condition (a) in the definition of $\sP _W (B)$, assertion 2 is trivially true from the definition of the lower density.

To see assertion 1, choose any $v \in T_{B,p}$ having norm 1, say with respect to the Bergman metric.  After a unitary change of coordinates in $\C^n$ (where the ball lies)  we may assume that that $v = c \frac{\di}{\di x^1}$, where $x^1,...,x^n$ are coordinates in $\C ^n$.
Consider the $(n-1,n-1)$-form
\[
\theta = \theta _v := (\ii)^{n-1} dx^2\wedge d\bar x ^2 \wedge \cdots \wedge dx^n \wedge d\bar x^n.
\]
We claim that $\theta \in \sP _W (B)$.  Indeed,  $\ii \di \dbar \theta = 0$ so condition (c) in the definition of $\sP _W(B)$ holds.  Condition (b) is clear in view of the local boundedness of $[W]_{\ve}$.  Condition (a) can be seen as follows:
\begin{eqnarray*}
\theta \wedge \omega _B &=& C \theta \wedge \left ( \frac{\omega _E}{1-|x|^2} + \frac{\ii \di |x|^2 \wedge \dbar |x|^2 }{(1-|x|^2)^2} \right )\\
&\ge & Ce^{-\lambda} \omega _E ^n \\
&=& C e^{n\lambda} \omega _B ^n.
\end{eqnarray*}

By the density condition there exists $\delta > 0$ and $r _o >> 0$
such that for all $r > r_o$,
\begin{eqnarray*}
1 - \delta & > & \frac{\left ( \Upsilon ^W _r + \frac{n}{n+1} \omega _B \right ) \wedge \theta}{\ii \di \dbar \kappa \wedge \theta}\\
&=&  \frac{\left ( \Upsilon ^W _r (v,\bar v)+ \frac{n}{n+1} \omega _B (v,\bar v) \right )}{\ii \di \dbar \kappa (v,\bar v)}.\\
\end{eqnarray*}
Observe that the density condition says this inequality holds uniformly on $B$.  Clearly if we rotate our original $v$ a little, this bound will still hold.  Since the unit sphere is compact, we can choose $r_o$ and $\delta$ so that the result holds for all $v$ in the unit sphere in $T_{B,p}$.  (Here, for the sake of simplifying the argument, we are exploiting the triviality of the tangent bundle $T_B$.)  This completes the proof.
\end{proof}

\subsection{A seemingly better notion of density}

In the paper \cite{osv}, a different notion of density was used.
The purpose of this section is to demonstrate the equivalence of
the density notions of the present paper and those in \cite{osv}.

Let us define the Bergman ball analogues of the densities used in \cite{osv}.  One first sets
\[
D_{z,r}(W,\kappa) := \sup _{v \neq 0} \frac{\Upsilon ^W _r (v,\bar v) + \frac{n}{n+1} \omega _B (v,\bar v)}{\ii \di \dbar \kappa (v,\bar v)}.
\]
Then one takes
\[
D^+(W,\kappa) := \limsup _{r \to 1} \sup _{z\in B} D_{z,r} (W,\kappa)
\]
and
$$
D^-(W,\kappa) := \liminf _{r \to 1} \inf _{z\in B} D_{z,r}
(W,\kappa).
$$

Note that $D_{z,r}(W,\kappa)$ is the maximum eigenvalue of
the $(1,1)$-form $\Upsilon ^W _r +\frac{n}{n+1} \omega _B$ with respect to the positive $(1,1)$-form $\ii \di \dbar \kappa$ at the point $z$.

\begin{lem}\label{form-rep}
Let $(M,\omega)$ be a Hermitian manifold of complex dimension $n$, and let $\alpha$ be a non-negative $(n-1,n-1)$ form on $M$.  Then for each $p$ there exists a vector $v \in T^{1,0}_{M,p}$ such that for
any real $(1,1)$-form $\beta$, one has
$$
\alpha \wedge \beta _p = \beta _p (v,\bar v) \omega ^n.
$$
\end{lem}

\noi The Lemma says that the mapping $\ii v \wedge \bar v \mapsto \theta _v$, with $\theta _v$ as in the proof of Lemma \ref{density-conditions}, is a pointwise isomorphism.

\begin{proof}
We shall use multi-linear algebra on $T_{M,p}$.  To this end, choose a unitary basis $e^1,...,e^n$ for $(T^* _{M,p})^{1,0}$ and $e_1,...,e_n$ its dual basis.  Let $\alpha ^{i\bar j}$ be a basis for $\Lambda ^{n-1,n-1} (T^*_{M,p})$ such that
$$
\ii e^k \wedge \bar e^{\ell} \wedge \alpha ^{i\bar j}= \delta ^{ik} \delta ^{\bar j \bar {\ell}} \frac{\omega ^n}{n!}.
$$
Let $A$ (resp. $B$) be the Hermitian matrix with entries $a_{i\bar j}$ (resp. $b_{i\bar j}$) such that at the point $p$,
$$
\alpha = a_{i\bar j}\alpha ^{i \bar j} \quad \left ( \text{resp.}\ \beta = b_{i \bar j} \ii e^i \wedge \bar e ^j \right ).
$$
Then
$$
\alpha \wedge \beta _p = \text{Trace}(AB^{\dagger}) \frac {\omega ^n}{n!} \quad \text{and} \quad \beta _p (v,\bar v) = v^{\dagger} B v.
$$
After a unitary rotation, we may assume that the basis $e^1,...,e^n$ diagonalizes $A$.  Thus, since $\alpha$ is positive, there exist non-negative numbers $\lambda _1,...,\lambda _n$ such that
$$
\text{Trace}(AB^{\dagger}) = \sum _{k=1} ^n \lambda _k b_{kk}.
$$
Taking
$$
v = \sum _{k=1} ^n \sqrt{\lambda _k} e_k
$$
completes the proof.
\end{proof}

\noi We can now obtain the following proposition.

\begin{prop}\label{density-relations} \
\begin{enumerate}
\item[(a)] $\sd ^+ (W,\kappa) = D^+(W,\kappa)$.
\item[(b)] $\sd ^- (W,\kappa) \le D^- (W,\kappa)$.
\end{enumerate}
\end{prop}

\begin{proof}
(a)  Fix $z \in B$ and $r \in [0,1)$.  By definition of $D^+(W,\kappa)$, we have that for any $\theta \in \sP _W (B)$,
\begin{eqnarray*}
D^+(W,\kappa) &\ge& D_{z,r} (W,\kappa) \\
&\ge &  \frac{\Upsilon ^W _r \wedge \theta (z) + \frac{n}{n+1} \omega _B \wedge \theta (z) }{\ii \di \dbar \kappa \wedge \theta (z)}.
\end{eqnarray*}
(In the second inequality we have used Lemma \ref{form-rep}.)
Taking the supremum over $z$ and then the $\limsup$ as $r \to 1$, we see that
$$
D^+(W,\kappa) \ge \sd ^+(W,\kappa) [\theta].
$$
Finally, taking the supremum of the right hand side over all $\theta \in \sP _W(B)$ shows that $D^+(W,\kappa) \ge \sd ^+ (W,\kappa)$.

To obtain the reverse inequality, fix $\ve > 0$.  Then for each $r <1$ sufficiently close to $1$ there exist $z \in B$ and $v \in \C ^n$ such that
\begin{eqnarray*}
D^+(W,\kappa) -\ve &\le &  \frac{\Upsilon ^W _r (v,\bar v) + \frac{n}{n+1} \omega _B (v,\bar v)}{\ii \di \dbar \kappa (v,\bar v)}\\
&=&  \frac{\Upsilon ^W _r \wedge \theta _v (z)  + \frac{n}{n+1} \omega _B \wedge \theta _v(z)}{\ii \di \dbar \kappa \wedge \theta _v(z)} +\ve \\
&\le & \sd ^+ (W,\kappa) [\theta _v] \le \sd ^+ (W,\kappa)  + \ve,
\end{eqnarray*}
where $\theta _v$ is defined as in the proof of Lemma \ref{density-conditions}.  The second-to-last inequality follows since $0 << r < 1$.  Since $\ve$ is arbitrary, 1 is proved.

\medskip

\noi (b)  Fix $\ve > 0$.  By definition of $\sd ^-(W,\kappa)$, there exists $\theta \in \sP _W(B)$ such that
$$
\sd ^- (W,\kappa) \le \sd ^- (W,\kappa)[\theta] + \frac{\ve}{2}.
$$
Moreover, by the definition of $\sd ^-(W,\kappa)[\theta]$ we have that for all $z\in B$ and all $r \in [0,1)$ sufficiently large,
$$
\sd ^- (W,\kappa)[\theta] \le \frac{\Upsilon _r ^W \wedge \theta (z) + \frac{n}{n+1} \omega _B \wedge \theta (z)}{\ii \di \dbar \kappa \wedge \theta (z)} + \frac{\ve}{2}.
$$
But by Lemma \ref{form-rep} and the definition of $D_{z,r} (W,\kappa)$,
$$
\frac{\Upsilon _r ^W \wedge \theta (z) + \frac{n}{n+1} \omega _B \wedge \theta (z)}{\ii \di \dbar \kappa \wedge \theta (z)} \le D_{z,r} (W,\kappa).
$$
This proves (b).
\end{proof}

\begin{thm}\label{sampling-relation}
$\sd ^- (W,\kappa) \ge D^-(W,\kappa).$
\end{thm}

\begin{proof}
We introduce the notation
$$
\Omega _{\delta} := \Upsilon ^W _r + \frac{n}{n+1} \omega _B - \left ( D^-(W,\kappa)-\delta \right ) \ii \di \dbar \kappa.
$$
Let $\delta > 0$ be given.  For $r >> 0$ we are going to construct a form $\theta \in \sP _W (B)$ such that
$$
\Omega _{\delta} \wedge \theta \ge 0.
$$
If this is done, the proof is complete.

By definition of $D^{-} (W,\kappa)$, there exists a locally finite open cover $U_j$ of $B$ and constant $(n-1,n-1)$-forms (i.e., forms of the type $\theta _v$ defined in the proof of Lemma \ref{density-conditions}) $\theta _j$ on $U_j$ such that
$$
\Omega _{\delta / 2} \wedge \theta _j \ge 0 \quad \text{on} \ U_j.
$$
By the uniform flatness of $W$ we may choose the cover $\{ U_j\}$ such that any point of $B$ is contained in some finite number of neighborhoods, this number depending only on the dimension.  Moreover, by the continuity of the forms $\Omega _{\delta}$ we may choose the forms $\theta _j$ so that if $U_j \cap U_k \neq \emptyset$ then $\theta _j - \theta _k$ is as small as we like.  In fact, by elementary anti-differentiation we may take forms $\mu _j$ depending quadratically on the (global) coordinates in $B$ such that $\theta _j = \ii \di \dbar \mu _j$ and  if $U_j \cap U_k \neq \emptyset$ then $||\mu _j - \mu _k||_{\sC ^2 (U_j \cap U_k)}$ is as small as we like, where $|| \cdot  ||_{\sC ^2}$ denotes $\sC ^2$-norm.

The argument we present here requires a little more precision.  Later we will have to control the size of the neighborhoods $U_j$ in order to make the $\theta _j-\theta _k$ small enough.  To this end, we choose the $U_j$ to be balls (or polydisks) of diameter $\ve$, measured with respect to the Bergman-Green distance.  We momentarily indicate this dependence on $\ve$ by writing $U_{j,\ve}$, $\mu _{j,\ve}$ and $\theta _{j,\ve}$.  Observe that if we take $\mu _{j,\ve}$ to be bihomogeneous quadratic in the Euclidean coordinates with origin that of $U_{j,\ve}$, then the uniform estimates for $\mu _{j,\ve}$ scale by $\ve ^2$, those for $D\mu _{j,\ve}$ by $\ve$, and those from $\theta _{j,\ve}$ are invariant with respect to $\ve$.

Let $\{ \psi _{j,\ve}\}$ be a partition of unity subordinate to the cover $\{ U_{j,\ve}\}$.  We may choose this partition so that
$$
\sum _j ||\psi _{j,\ve} \mu _{j,\ve} ||_{\sC ^2} \le C
$$
for some constant $C$ independent of $\ve$.  Indeed, as the neighborhoods $U_{j,\ve}$ scale by $\ve$, the estimates for $D\psi _{j,\ve}$ scale by $\ve ^{-1}$ while those for $D^2 \psi _{j,\ve}$ scale by $\ve ^{-2}$.  Thus the desired estimate follows from the product rule
$$
D^2(\psi _{j,\ve} \mu_{j,\ve} ) = \mu _{j,\ve } D^2 \psi _{j,\ve} + (D\mu _{j,\ve})(D\psi _{j,\ve} ) + \psi _{j,\ve} D^2 \mu _{j,\ve}.
$$
Thus is is clear that we have scale invariant estimates.  To simplify the notation, we can now drop the notational dependence on $\ve$.

We would like to correct the local forms $\theta _j$ so that they can be pieced together to give us an element of $\sP _W(B)$ with the desired density.  We shall use cocycles to do this.  To this end, the obstruction to the $\theta _j$ piecing together to give a global form is carried by the 1-cocycle
$$
\alpha _{jk} = \theta _j - \theta _k = \ii \di \dbar (\mu _j - \mu _k)
$$
supported on $U_j \cap U_k$.  By our choice of the $\theta _j$,  the $\alpha _{jk}$ are small in $\sC ^0$-norm.  We now define
$$
\eta _j = \ii \di \dbar \left ( \sum _k \psi _k (\mu _{j}- \mu _{k}) |U_j \cap U_k \right ).
$$
By modifying our choices of the $\mu _j$ we may make the $\eta _j$ as small as we like.  Moreover,
$\ii \di \dbar \eta _j = 0$ and
\begin{eqnarray*}
\eta _j - \eta _\ell &=& \ii \di \dbar \sum _k \psi _k (\mu _j - \mu _k + \mu _k - \mu _{\ell}) \\
&=& \ii \di \dbar \sum _k \psi _k ( \mu _j - \mu _{\ell}) \\
&=& \alpha _{j \ell}.
\end{eqnarray*}
It follows that
$$
\theta = \theta _j - \eta _j \quad \text{on} \quad U_j
$$
is well defined and belongs to $\sP _W (B)$.  Moreover, by choosing the $\mu _j-\mu _k$ even smaller if necessary, we see that
$$
\Omega _{\delta} \wedge \theta \ge 0,
$$
as desired.
\end{proof}

\section{Interpolation}\label{interp-section}

\subsection{A negative function singular along a hypersurface}

Recall that
\[
V_n(r) := \int _{B(0,r)} \omega _B^n.
\]
As the Bergman metric is invariant under automorphisms, one sees that for each $a \in B$, $V_n(r)$ is also the Bergman volume of Bergman-Green balls $E(a,r) = F_a(B(0,r))$.

Let
\begin{eqnarray*}
\Gamma_r(z,\zeta)=G_B(z,\zeta) - \frac{1}{V_n(r)} \int_{E(z,r)}
G_B(x,\zeta)\omega_B^n(x).
\end{eqnarray*}
Since $\Delta _B G(\cdot ,\zeta) \equiv 0$ on $B - \{ \zeta \}$,
we see from Corollary \ref{berg-ineq} that $\Gamma_r$ is
non-negative and is supported on the set
\[
\left \{(z,\zeta) \in B \times B \ \left |\ |F_z(\zeta)| \le r \right
. \right \},
\]
which contains a neighborhood of the diagonal in $B \times B$.

We define the function
\begin{eqnarray*}
s_r(z)&:=& \int_B \Gamma _r (z,\zeta)\omega_B^{n-1}(\zeta)
\wedge dd^c \log|T|^2 (\zeta)\\
&=& \int_{\{ \zeta \ ;\ |F_z(\zeta)|<r\}} \Gamma _r (z,\zeta) \omega_B^{n-1}(\zeta) \wedge dd^c\log|T|^2 (\zeta).
\end{eqnarray*}
By the Lelong-Poincar\'e identity,
\begin{equation}\label{pl-id}
s_r(z)=2\pi  \int_{W_{z,r}} \left ( G_B (z,\zeta) - \frac{1}{V(r)}
\int_{E(z,r)} G_B(x,\zeta)\omega_B^n(x) \right) \omega_B^{n-1}(\zeta),
\end{equation}
where
\begin{displaymath}
W_{z,r}=W \cap \left \{ \zeta \ ; \ |F_z(\zeta)|<r \right \}
\end{displaymath}

\begin{prop}\label{log-sing}
Let $T \in \co(B)$ be a holomorphic function so that $W=T^{-1}(0)$ and $dT$ is nowhere zero on $W$. Then
\[
s_r(z)=\log|T(z)|^2-\frac{1}{V_n(r)}\int_{E(z,r)} \log |T(\zeta)| ^2 \omega_B^{n}(\zeta).
\]
In patricular,
\[
\frac{1}{2\pi} dd^c s_r (z) = [W] - \Upsilon ^W _r (z).
\]
\end{prop}

\begin{proof}
Let $\vp \in \sC ^{\infty} _0 (B)$ be a function whose total integral with respect to Euclidean volume is $1$, and let $\chi _{\ve}$ be the characteristic function of the set $\{ z\in B \ ;\ |z| < 1-2\ve \}$.  Let $\vp _{\ve} (x) = \ve ^{-2n}\vp (\ve ^{-1} x)$, and set
\[
f_{\ve} = (\chi _{\ve} \log |T|^2 ) * \vp _{\ve}.
\]
Then $f_{\ve}$ is smooth with compact support in $B$, and
\[
s_r (z) = \lim _{\ve \to 0} \int_B \Gamma _r (z,\zeta)\omega_B^{n-1}(\zeta) \wedge dd^c f_{\ve} (\zeta).
\]
But by definition of Green's function,
\[
\int _B G(z,\zeta) \omega_B^{n-1}(\zeta) \wedge dd^c f_{\ve} (\zeta) = \int _B G(z,\zeta) dd^c ( f_{\ve} \omega_B^{n-1}) (\zeta) ) = f_{\ve} (z).
\]
The proof is completed by letting $\ve \to 0$.
\end{proof}

\begin{lem}\label{s-prop}
The function $s_r(z)$ has the following properties:
\begin{enumerate}
\item[1.] It is non-positive.

\item[2.] For each $r,\epsilon >0$ there exist a constant
$C_{r,\epsilon}$ such that if $\delta _B(z,W)>\epsilon$ then
$s_r(z)>-C_{r,\epsilon}$.

\item[3.] The function $e^{-s_r}$ is not locally integrable at any point of W.
\end{enumerate}
\end{lem}

\noi Here $\delta _B (z,W) = \inf \{ |F_z(w)|\ ;\ w\in W\}$.

\begin{proof}
By Corollary \ref{berg-ineq} and the fact that $\Delta _B G ( \cdot , \zeta ) \equiv 0$ on $B - \{ \zeta \}$, $\Gamma_r \leq 0$
and 1 follows.  Moreover, 3 is an immediate consequence of
Proposition \ref{log-sing}.

To see 2, we first note that since $\delta _B (z,W) > \ve$,  $G_B(z,\zeta) > A_{\ve}$.  Thus it suffices to obtain an estimate
$$
-\int_{E(z,r)} G_B (x,\zeta) \omega_B^n (x) = -\int _{B(0,r)} G_B(x,y) \omega ^n_B(x) \le D_r
$$
for some $D_r >0$ and all $y = F_z(\zeta) \in B(0,r)$.  To do this, it is enough to estimate the integral
$$
I(r) := -\int _{B(0,(r+1)/2)} G_B(x,y) \omega ^n_B(x).
$$

Fix $y \in B(0,r)$.  Let $\rho>0$ be the largest number such that
$$
B(y,\rho) \subset B(0,(r+1)/2).
$$
One has
$$
n_r \le \rho \le \frac{r+1}{2}
$$
for some $n_r >0$ depending on $r$ but not on $y$.

Write
$$
I(r) = I_1(r) +I_2(r),
$$
where
$$
I_1(r) := -\int _{E(y,\rho)} G_B(x,y) \omega ^n_B(x)
$$
and
$$
I_2(r) := -\int _{B(0,(r+1)/2) - E(y,\rho)} G_B(x,y) \omega ^n_B(x).
$$
Now
$$
I_1(r) = - \int _{B(0,\rho)} \gamma _B (x) \omega ^n _B (x)
$$
is clearly bounded by a constant independent of $y$.

Next, note that for
$$
x \in B(0,(r+1)/2) - E(y,\rho)
$$
one has the estimate
$$
|F_y(x)| \ge \rho \ge n_r.
$$
It follows that for such $x$,
$$
G(x,y) \ge - N_r
$$
for some $N_r \in \R$ independent of $y$.  Thus
\begin{eqnarray*}
I_2(r) &\le&  N_r \int _{B(0,(r+1)/2) - E(y,\rho)} \omega^n _B(x)\\
&\le & N_r \int _{B(0,(r+1)/2)} \omega ^n_B(x),
\end{eqnarray*}
and the latter is independent of $y$.  Thus 2 follows.
\end{proof}

\begin{rmk}
There is a direct proof of Lemma \ref{s-prop}.3 that does not use the formula of Proposition \ref{log-sing}.  Since we will make use of the calculation needed, we present this proof now.

We may assume that $W$ is the coordinate hyperplane $z_n =0$ and
$z = z^n e_n$ for $|z^n | \le \ve$ with $\ve$ sufficiently small.
(Though we do not use it here, later we will exploit the fact
that, by the uniform flatness of $W$, $\ve > 0$ may be taken
independent of the point on $W$ which has been translated to the
origin.)  Let $U$ be a sufficiently small neighborhood of the
origin.  Using the formula \eqref{pl-id}, we estimate that
\[
s_r(z) = 2\pi \int _{U\cap W}\Gamma _r(z,\zeta) \omega _B ^{n-1}(\zeta) + O(1).
\]
The same method used in the proof of Lemma \ref{s-prop}.2 allows us to estimate the part of $\Gamma _r$ involving the integral, so we may replace $\Gamma _r$ be the Green's function.  Letting $\omega _B ^{n-1} = 2 r^{2n-3} dr d\sigma _{2n-3}$ be the decomposition into polar coordinates and setting
\[
A_{n-1} = \int _{S^{2n-3}} d\sigma _{2n-3} = (2\pi) ^{n-1},
\]
we obtain from the form of the Green's function that
\[
s_r(z) = 2\pi C_n (n+1)^{n-1} A_{n-1} \int _0 ^{\alpha}-\frac {2r^{2n-3}dr}{\left ( r^2+ |z^n|^2 \right ) ^{n-1}} + O(1) =\log |z^n |^2 +O(1)
\]
where $\alpha > 0$ is a sufficiently small number depending on $\ve$.
\qed
\end{rmk}

\subsection{The proof of Theorem \ref{interp-thm}} \label{ot-pf}

We fix a compact subset $\Omega \relcomp B$.  This set will be fixed until the last part of the argument, when we let $\Omega \to B$.

Let
\[
\sigma _r = s_r  - \lambda  - \sup _{\Omega} (-\lambda).
\]
Note that $\sigma _r \le 0$.

\subsubsection*{Tubular limits}

For each $\Omega \relcomp B$, let
\[
\Omega _{\ve} := \Omega \cap \{ \sigma _r < \log \ve ^2 \}.
\]
\begin{lem}\label{limits}
Let $W \subset B$ be uniformly flat.  Then there exists a positive constant $C>0$ such that for all $\Omega \relcomp B$ and all $f$ holomorphic in a neighborhood of $\overline \Omega$,
\[
\limsup _{\ve \to 0} \frac{1}{\ve^2} \int _{\Omega _{\ve}} |f|^2 e^{-\kappa} \omega _B ^n \le C \int _{\Omega \cap W} |f|^2 e^{-\kappa} \omega _B ^{n-1}.
\]
\end{lem}

\begin{proof}[Sketch of proof]
We may assume the right hand side is finite. Moreover, we can take $\Omega = E(a,\delta)$ for some $a \in W$, with $\delta$ so small that $W \cap \Omega$ is the graph of a quadratic hypersurface.  By uniform flatness, $\delta$ can be taken independent of $a$.

Consider first the case $a=0$.  Then $\Omega =B(0,\delta)$, and the result follows after an elementary analysis of the properties of $s_r$ as in the proof of Lemma \ref{s-prop}, and the remark following that proof.

If we now apply the automorphism $F_a$ to $B(0,\delta)$, then Lemma \ref{u-f-char} and the Aut$(B)$-invariance of $\omega _B$ show that the same estimates hold on $E(a,\delta)$.
\end{proof}

\subsubsection*{The twisted Bochner-Kodaira Technique}

We fix a smoothly bounded pseudoconvex domain $\Omega \relcomp B$.  Let us denote by $\dbar ^* _{\nu}$ the formal adjoint of $\dbar$ in the Hilbert space of $(0,1)$-forms on $\Omega$, square integrable with respect to a weight $e^{-\nu} \omega _B ^n$. For a $(0,1)$-form $u = u_{\bar \alpha} d\bar z ^{\alpha}$, one has
$$
\dbar ^* _{\nu} u = - e^{\nu + (n+1)\lambda} \di _{\alpha}\left ( e^{-(\nu +(n+1)\lambda)} u^{\alpha} \right ).
$$
Recall that for $(0,1)$-forms $u$ in the domains of $\dbar$ and $\dbar ^* _{\nu}$, Bochner-Kodaira Identity is
\begin{eqnarray}
\nonumber && \int _{\Omega} \left | \dbar ^* _{\nu} u \right | ^2 e^{-\nu} \omega _B ^n + \int _{\Omega} \left | \dbar u \right | ^2
e^{-\nu} \omega _B ^n\\
\label{bk} &&= \int _{\Omega}  \left ( \left (\di _{\alpha}\di
_{\bar \beta}(\nu +(n+1)\lambda)\right ) u^{\alpha} \overline
{u^{\beta}}\right ) e^{-\nu} \omega _B ^n + \int _{\Omega} \left | \overline \nabla u \right | ^2 e^{-\nu} \omega _B ^n \\
\nonumber &&\qquad \qquad  + \int _{\di \Omega} \left ( \di
_{\alpha} \di _{\bar \beta} \rho \right ) u^{\alpha} \overline{u^{\beta}} e^{-\nu} d^{c}(-(n+1)\lambda ) \wedge \omega _B ^{n-1},
\end{eqnarray}
where $\rho$ is a defining function for $\Omega$ such that
$|d\rho|\equiv 1$ on $\di \Omega$.  (See, for example,
\cite{siu-82}.)  The term $(n+1)\lambda$ in the first integral on the right hand side of \eqref{bk} comes from the Ricci
curvature of $\omega _B$.  Writing
$$
e^{-\psi} = \frac{e^{-\nu}}{\tau}
$$
we obtain
$$
\dbar ^* _{\nu } u = \dbar ^* _{\psi} u - \frac{\left ( \di
_{\alpha} \tau \right ) u^{\alpha} }{\tau} \quad {\rm and} \quad \di _{\alpha} \di _{\bar \beta} \psi = \di _{\alpha} \di _{\bar \beta} \nu + \frac{\di _{\alpha} \di _{\bar \beta} \tau }{\tau} - \frac{\left ( \di _{\alpha}\tau \right ) \left (  \di _{\bar
\beta} \tau \right )}{\tau ^2}.
$$
Substitution into (\ref{bk}), followed by some simple
manipulation, gives the

\medskip

\noi {\sc Twisted Bochner-Kodaira Identity} for $(0,1)$-forms:  If $u$ is a $(0,1)$-form in the domain of $\dbar ^*$, then
\begin{eqnarray}\label{tbk}
&& \int _{\Omega} \tau \left | \dbar ^* _{\psi} u \right |^2
e^{-\psi} \omega _B ^n + \int _{\Omega} \tau \left | \dbar u
\right |^2 e^{-\psi} \omega _B ^n\\
\nonumber &=& \int _{\Omega}  \left ( \tau \left ( \di
_{\alpha}\di _{\bar \beta} (\psi +(n+1)\lambda ) \right)
u^{\alpha} \overline {u^{\beta}} - \left ( u^{\alpha} \overline
{u^{\beta}}\di _{\alpha}\di _{\bar \beta} \tau  \right ) \right .\\
\nonumber && \qquad  \qquad  + \left . 2 \re \left ( \left ( \di
_{\alpha} \tau \right ) u^{\alpha} \overline{ \dbar ^* _{\psi} u } \right ) \right ) e^{-\psi} \omega _B ^n + \int _{\Omega} \tau
\left | \overline \nabla u \right | ^2 e^{-\psi} \omega _B ^n \\
\nonumber && \qquad + \int _{\di \Omega} \tau \left ( \di
_{\alpha} \di _{\bar \beta} \rho \right ) u^{\alpha}
\overline{u^{\beta}} e^{-\psi} d^{c}(-(n+1) \lambda) \wedge \omega _B ^{n-1}.
\end{eqnarray}
We now use positivity of the last two integrals on the right hand side, together with the Cauchy-Schwarz inequality applied to the first term in the third line, to obtain the so-called

\medskip

\noi {\sc Twisted basic estimate:} If $u$ is a $(0,1)$-form in the
domain of $\dbar ^*$, then
\begin{eqnarray}\label{tbe}
&& \int _{\Omega} \left ( \tau + A\right ) \left | \dbar ^* _{\psi} u \right |^2 e^{-\psi} \omega _B ^n + \int _{\Omega} \tau \left | \dbar u \right |^2 e^{-\psi} \omega _B ^n\\
\nonumber && \qquad \ge  \int _{\Omega}  \left ( \tau \left (\di
_{\alpha}\di _{\bar \beta}( \psi +(n+1)\lambda )  \right )u ^{\alpha} \overline {u^{\beta}} \right . \\
\nonumber && \qquad \qquad \quad  \left . - \di _{\alpha} \di _{\bar \beta} \tau u ^{\alpha} \overline{u^{\beta}} - \frac {1}{A} \left | \left ( \di _{\alpha} \tau \right ) u^{\alpha} \right | ^2 \right ) e^{-\psi} \omega _B ^n .
\end{eqnarray}

\subsubsection*{Choice of $\psi$, $\tau$ and $A$}
From the very beginning, we choose
\[
\psi = \kappa + \sigma _r.
\]
By the density hypothesis (via Lemma \ref{density-conditions}.1) and the fact that $\ii \di \dbar s_r = [W] - \Upsilon ^W _r$, one has
\begin{eqnarray*}
\ii \di \dbar (\psi + (n+1) \lambda ) &=& \ii \di \dbar (\kappa +
n \lambda +s_r) \\
& \ge & c \ii \di \dbar \kappa.
\end{eqnarray*}

Next, fix $\gamma > 1$.  We define
$$
\xi = \log \left ( e^{\sigma _r} +\ve ^2 \right ),
$$
with $\ve >0$ so small that $\gamma - \xi \ge 1$.  One has
\begin{eqnarray*}
&& \ii \di \dbar \xi \\
&=& \ii \di \left ( \frac{e^{\sigma _r}}{e^{\sigma _r} + \ve ^2 } \dbar \sigma _r \right ) \\
&=& \frac{e^{\sigma _r}}{e^{\sigma _r} + \ve ^2 } \ii \di \dbar \sigma _r + \frac{\ve ^2}{(e^{\sigma _r} + \ve ^2 )^2} e^{\sigma _r} \ii \di \sigma _r \wedge \dbar \sigma _r \\
&=& \frac{e^{\sigma _r}}{e^{\sigma _r} + \ve^2}\left ( \tfrac{1}{n+1} \omega _B -  \Upsilon ^W _r \right ) + \frac{4 \ve ^2}{(e^{\sigma _r} + \ve ^2 )^2} \ii \di \left ( e^{\frac{1}{2} \sigma _r} \right ) \wedge \dbar \left ( e^{\frac{1}{2} \sigma _r} \right ),
\end{eqnarray*}
where the last equality follows since $\ii \di \dbar \sigma _r = [W] +\tfrac{1}{n+1} \omega _B - \Upsilon ^W _r$ and $e^{\sigma _r}|W \equiv 0$.

Let $0<\alpha << 1$ and set
$$
a = \gamma - \alpha \xi.
$$
Observe that $a \ge 1$.  Moreover, we have
\begin{eqnarray*}
 &&- \ii \di \dbar a \\
&=& \alpha \ii \di \dbar \xi \\
&=& \frac{\alpha e^{\sigma _r}}{e^{\sigma _r} + \ve^2}\left (
\tfrac{1}{n+1} \omega _B -  \Upsilon ^W _r \right ) + \frac{4
\alpha \ve^2}{(e^{\sigma _r} + \ve ^2)^2} \ii \di \left
(e^{\tfrac{1}{2}\sigma _r} \right ) \wedge \dbar \left
(e^{\tfrac{1}{2} \sigma _r} \right ).
\end{eqnarray*}

Now let
\[
\tau = a+ \log a \quad \text{and} \quad A= (1+a)^2.
\]
Then $\tau \ge 1$ and we have
\[
\di \tau = \left ( 1 + \frac{1}{a} \right ) \di a  \quad \text{and} \quad \ii \di \dbar \tau = \left ( 1+ \frac{1}{a} \right ) \ii \di \dbar a - \frac{1}{a^2} \ii \di a \wedge \dbar a,
\]
and thus
$$
- \ii \di \dbar \tau - \frac{\ii \di \tau \wedge \dbar \tau }{A} = \left ( 1+\frac{1}{a}\right ) \left ( - \ii \di \dbar a\right ) \ge - \ii \di \dbar a.
$$
It follows that
\begin{eqnarray*}
&& \tau \ii \di \dbar (\psi + (n+1)\lambda) - \ii \di \dbar \tau - \frac{|\di \tau |^2}{A} \\
&& \ge c \ii \di \dbar \kappa + \frac{\alpha e^{\sigma _r}}{e^{\sigma _r} + \ve^2}\left ( \tfrac{1}{n+1} \omega _B -  \Upsilon ^W _r \right ) \\
&& \qquad  + \frac{4\alpha \ve ^2}{(e^{\sigma _r} + \ve ^2 )^2} \ii\di \left ( e^{\frac{1}{2} \sigma _r}\right ) \wedge \dbar \left ( e^{\frac{1}{2} \sigma _r}\right )\\
&& \ge  \frac{4\alpha \ve ^2}{(e^{\sigma _r} + \ve ^2 )^2} \ii\di \left ( e^{\frac{1}{2}\sigma _r}\right ) \wedge \dbar \left ( e^{\frac{1}{2} \sigma _r}\right ),\\
\end{eqnarray*}
provided we take $\alpha$ sufficiently small.  (For example, by the density hypothesis as rephrased in Lemma \ref{density-conditions}.1 we may take $\alpha = c$.)  Substituting into the twisted basic estimate (\ref{tbe}), we obtain the following lemma.

\begin{lem}\label{hormander}
If $u$ is a $(0,1)$-form in the domain of $\dbar ^*$, then
\begin{eqnarray*}
\int _{\Omega} \left (\tau + A \right ) \left | \dbar ^* _{\psi} u \right | ^2 e^{-\psi} \omega _B ^n + \int _{\Omega} \tau \left | \dbar u \right | ^2 e^{-\psi} \omega _B ^n \\
\qquad \ge c  \int _{\Omega} \frac{4\ve ^2}{(e^{\sigma _r} + \ve
^2)^2} \left | \di \left ( e^{\frac{1}{2}\sigma _r} \right )(u) \right
|^2 e^{-\psi} \omega _B ^n.
\end{eqnarray*}
\end{lem}

\subsubsection*{An a priori estimate}\label{a-priori-section}

We write $\Omega _j = B\left (0,\frac{j}{1+j}\right )$.  Suppose
we are given $f \in \fH ^2 (W,\kappa)$.  Since $W$ is a closed
submanifold of $B$, there exists a holomorphic extension $\tilde
f$ of $f$ to $B$. We write
$$
W_j = W \cap \Omega _j, \quad f_j = f|W_j \quad {\rm and} \quad \tilde f _j = \tilde f |\Omega _j.
$$
Observe that
$$
\int _{W_j} |f_j|^2 e^{-\kappa} \omega _B ^{n-1} \le \int _W |f|^2 e^{-\kappa } \omega _B ^{n-1} < +\infty.
$$
Let $\chi \in \cc ^{\infty} _0 ([0,1))$ be such that
$$
0 \le \chi \le 1, \quad  \chi |[0,1/3] \equiv 1 \quad {\rm and}
\quad \sup _{[0,1)} |\chi '|\le 2.
$$
We set
\[
\chi _{\ve} = \chi \left ( \frac{e^{\sigma _r}}{\ve ^2} \right )
\]
and define the 1-forms $\alpha _{\ve,j}$ on $\Omega _j$ by
\[
\alpha_{\ve ,j} = \dbar \chi _{\ve} \tilde f_j.
\]
We note that for $\ve$ sufficiently small, $\alpha _{\ve ,j}$ is
supported on the tubular neighborhood
\[
\Omega _{\ve ,j} := \Omega _j \cap \left \{ e^{\frac{1}{2} \sigma _r} \le \ve \right \}
\]
of $W_j$ in $\Omega _j$.  Thus, for a $(0,1)$-form $u$ with
compact support on $\Omega _j$, we have
\begin{eqnarray*}
\left | ( \alpha _{\ve ,j} , u) \right |^2 &=& \left | \int _{\Omega _j} \left < \chi ' \left ( \frac{e^{\sigma _r}}{\ve ^2} \right ) \frac{\dbar (e^{\sigma _r})}{\ve ^2} \tilde f_j , u \right > e^{-\psi} \omega _B ^n \right | ^2\\
&\le & \left ( \frac{2}{\ve ^2} \int _{\Omega _j} \left |\chi ' \left ( \frac{e^{\sigma _r}}{\ve ^2} \right ) \right | \left | \di \left ( e^{\frac{1}{2} \sigma _r} \right ) (u) \right | |\tilde f_j| e^{-\frac{1}{2} \sigma _r} e^{-\kappa }\omega _B ^n \right ) ^2\\
& \le & \frac{16}{\ve ^4} \left ( \int _{\Omega _{\ve,j}} |\tilde
f _j |^2 \frac{(e^{\sigma _r} +\ve ^2 )^2}{4\ve ^2} e^{-\kappa }
\omega _B ^n \right ) \\
&& \qquad  \times \int _{\Omega _j} \frac{4\ve ^2}{(e^{\sigma _r}
+ \ve ^2)^2} \left | \di \left ( e^{\frac{1}{2}\sigma _r} \right ) (u) \right |^2 e^{-\psi} \omega _B ^n \\
&\le & \frac{16}{c} C_{\ve,j} \left ( || T^* u||^2 + ||Su||^2
\right ) ,
\end{eqnarray*}
where
$$
Tu = \dbar \left ( \sqrt{\tau + A} u \right ) \quad {\rm and}
\quad Su = \sqrt \tau \left ( \dbar u \right ),
$$
and
$$
C_{\ve ,j} = \frac{1}{\ve ^2} \int _{\Omega _{\ve ,j}} |\tilde f_j |^2  e^{-\kappa} \omega _B^{n}.
$$
Thus the last inequality follows from Lemma \ref{hormander} and the fact that $e^{\sigma _r} < \ve ^2$ on $\Omega _{\ve, j}$.

\medskip

By standard Hilbert space methods, we have the following $L^2$
twisted-$\dbar$ theorem.

\begin{thm}\label{T-soln}
There exists a function $h_{j,\ve}$ on $\Omega _j$ such that
$$
T h_{\ve ,j} = \alpha _{\ve ,j} \quad {\rm and} \quad \int
_{\Omega _j} |h_{\ve ,j} |^2 e^{-\psi} \omega _B^n \le
\frac{16}{c} C_{\ve,j}.
$$
In particular, $h_{\ve ,j}|W \equiv 0.$
\end{thm}

\begin{proof}
Consider the linear functional
\[
\sL : T^* u \mapsto (u, \alpha_{\ve , j}),
\]
where $u \in \text{Kernel} S \cap {\rm Domain} (T^*)$.  The estimate
\[
|(u, \alpha_{\ve , j})|^2 \le \frac{16C_{\ve,j}}{c} ||T^*u||^2
\]
($Su = 0$) means $\sL$ is continuous on the image of $T^*$, hence on the closure of that image.  Extend $\sL$ to be $0$ in Image$(T^*) ^{\perp}$.  Then $\sL$ is a continuous linear functional in our Hilbert space, and thus, by the Riesz Representation Theorem, is represented by some element $h_{\ve ,j}$ having the same norm as $\sL$   in the orthogonal direction.  Elliptic regularity implies that $h_{\ve , j}$ is smooth.

It remains only to prove the assertion about the vanishing of $h_{\ve ,j}$.  But by Lemma \ref{s-prop}.3, $e^{-\psi}$ is not locally integrable at any point of $W$, and thus the vanishing of $h_{\ve ,j}|W$ follows.
\end{proof}

\subsubsection*{Conclusion of the proof of Theorem \ref{interp-thm}}
Observe first that by Lemma \ref{limits} there exists a constant $C>0$ such that, for all $j$,
\[
\limsup _{\ve \to 0} C_{\ve,j} \le C \int _W |f|^2 e^{-\kappa}
\omega _B ^{n-1}.
\]
We set
\[
F_{\ve,j} = \chi _{\ve} \tilde f _j -  \sqrt{(\tau +A)} \ h_{\ve
,j} \quad {\rm on} \ \Omega _j.
\]
By Theorem \ref{T-soln}, $F_{\ve ,j}$ is holomorphic on $\Omega _j$ and $F_{\ve ,j}|W_j - f_j \equiv 0$.  Moreover there exists a constant $M$ such that
$$
\int _{\Omega _j} |F_{\ve ,j} |^2 e^{-\kappa} \omega _B ^n \le
M\left ( o(1) + \int _W |f|^2 e^{-\kappa } \omega _B ^{n-1} \right ), \quad \ve \sim 0.
$$
Indeed, the integral
$$
\int _{\Omega _j} |\chi _{\ve} \tilde f_j |^2 e^{-\kappa}\omega _B^n
$$
is negligible for small $\ve$, since the integrand is locally integrable and supported on a set of arbitrarily small measure. On the other hand,
\begin{eqnarray*}
\int _{\Omega _j} (\tau +A) |h_{\ve,j}|^2 e^{-\kappa}\omega _B ^n &=& \int _{\Omega _j} e^{\sigma _r} (\tau +A)
|h_{\ve,j}|^2 e^{-\psi}\omega _B ^n \\
&\le & \left ( \sup _{\Omega_j } e^{\sigma _r} (\tau + A) \right )
\int _{\Omega _j} |h_{\ve ,j} |^2 e^{-\psi} \omega _B ^n \\
&\le & C_{\ve ,j} e^{\frac{\gamma}{\alpha }} \left ( \sup
_{\Omega_j } e^{- \frac{1}{\alpha } a} (a + \log (a)+(1+a)^2) \right ) \\
&\le & K C_{\ve ,j}
\end{eqnarray*}
for some universal constant $K$ depending only on the density of $W$.  The last estimate holds since $a \ge 1$.

By Corollary \ref{berg-ineq} and the Lebesgue Dominated Convergence Theorem,
\[
F_j = \lim _{\ve \to 0} F_{\ve , j}
\]
exists, uniformly for each fixed $j$. Moreover, since $F_{\ve ,j}
= f$ on $W_j$ and $F_{\ve ,j} \to F_{j}$ pointwise, we have $F_{j}=f$ on $W_j$ for all $j$. We thus have a sequence of functions $F_j$,
holomorphic by Montel's Theorem, such that $F_j |W =f$ and
\[
\int _{\Omega _j} |F_j|^2 e^{-\vp} \omega ^n \le C \int _W |f|^2 e^{-\vp} \omega ^{n-1}.
\]
Moreover, the constant $C$ does not depend on $j$.  Letting $j \to \infty$, we obtain (again by corollary \ref{berg-ineq}, the Dominated Convergence Theorem and Montel's Theorem) a holomorphic function $F$ that also agrees with $f$ on $W$, and furthermore satisfies
\[
\int _B |F|^2 e^{-\kappa} \omega _B ^n \le C \int _W |f|^2
e^{-\kappa}.
\]
This completes the proof of Theorem \ref{interp-thm}

\section{Sampling}\label{samp-section}

\subsection{A construction of Berndtsson-Ortega Cerd\` a.}

\begin{lem}\label{quimbo-trick}
Let $\vp$ be a subharmonic function on the unit disk $\D$.  Then there exist a positive constant $K$ and a holomorphic function on $G \in \co (\D (0,1/2))$ such that $G(0)=0$ and
\[
\sup _{\D \left (0,1/2\right )} \left |\vp - \vp (0) - 2\re \ G \right | \le K.
\]
Moreover, if $\vp$ depends smoothly on a parameter, then so does $G$.
\end{lem}

\noi The proof of this lemma, which uses Riesz Potentials, can be found in \cite{quimbo}.

\subsection{Restriction from tubes and the upper inequality}

\begin{prop}\label{tube-smpl}
Let $W$ be a uniformly flat smooth hypersurface.  Then there
exists a constant $C>0$ such that for all $\ve >0$ sufficiently
small and all $F \in \sh ^2(N_{\ve}^B(W),\kappa)$ one has
$$
C \ve^2 \int_W |F|^2 e^{-\kappa} \omega_B^{n-1} \leq \int_ {N_{\ve}^B(W)} |F|^2 e^{-\kappa} \omega_B^n.
$$
\end{prop}

\begin{proof}
Let
\[
D(0,\ve) =  \{ (z',z^n)\in \C ^{n-1} \times \C \ ;\ |z'| \le \ve, \ |z^n| < \ve \}.
\]
Via Lemma \ref{u-f-char}, the uniform flatness of $W$ implies that $N^W _{\ve} (W)$ is a union of open sets $U_j$ such that for each $j$ there is some $F_{z_j} \in Aut (B)$ for which
\[
F_{z_j} (U_j) \sim D(0,\ve).
\]
Moreover, this approximation may be taken uniform in $j$.  Thus it suffices to prove that for some $\ve > 0$ and all $a \in W$,
\[
C \ve ^2 \int _{D(0,\ve)\cap F_a(W)} |F |^2 e^{-\kappa} \omega _B ^{n-1} \le \int _{D(0,\ve)} |F|^2 e^{-\kappa} \omega _B ^n.
\]
After a change of variables provided by Lemma \ref{u-f-char}, we may assume that $F_a(W) \subset \C ^{n-1} \times \{0\}$.

Now, by Lemma \ref{quimbo-trick} there exists a function $G$, holomorphic in $z^n$, such that
\[
G(z',0)| \equiv 0 \quad {\rm and} \quad e^{-\kappa (z',0) + 2\re
G(z',z^n)} \le c e^{-\kappa (z',z^n)}
\]
for some $c>0$. We then have
\begin{eqnarray*}
\int _{B(0,\ve) \cap \C ^{n-1} \times \{0\}} |F |^2
e^{-\kappa}\omega _B ^{n-1} &=& \int _{B(0,\ve) \cap \C ^{n-1} \times \{0\}} |F e^{G} |^2 e^{-\kappa}\omega _B ^{n-1} \\
& \le & \frac{C_o}{\ve ^2} \int _{D(0,\ve)} |F
e^{G} |^2 e^{-\kappa(z',0)}\omega _B ^n \\
&=& \frac{C_o}{\ve ^2} \int _{D(0,\ve)} |F |^2
e^{-\kappa(z',0) + 2 \re \ {G}}\omega _B ^n \\
&\le & \frac{1}{C \ve^2}\int _{D(0,\ve)} |F|^2 e^{-\kappa (z',z^n)}\omega _B ^n. \\
\end{eqnarray*}
The first inequality follows from the sub-mean value property for radial measures in the disk (see also Corollary \ref{berg-ineq}).  This completes the proof.
\end{proof}

\begin{cor}\label{upper-sampling-ineq}
If $W$ is a uniformly flat hypersurface then there exists a
constant $M$ such that for all $F \in H^2 (B,\kappa)$,
$$
\int _W |F |^2 e^{-\kappa} \omega _B ^{n-1} \le M \int _B |F |^2 e^{-\kappa} \omega _B ^n.
$$
\end{cor}

\subsection{Regularization of the singular function $s_r$}

Consider the function
$$
s_{r,\ve}(z) := \frac{1}{V_n(\ve)} \int _{E(z,\ve)} s_r \omega _B
^n.
$$
In this section we prove the following result.

\begin{lem}\label{curb}
The function $s_{r,\ve}$ enjoys the following properties.
\begin{enumerate}
\item
$$
\lim _{\ve \to 0} dd^c s_{r,\ve} = [W] - \Upsilon ^W_r.
$$

\item For each $r$ there exists a constant $C_r$ such that if $0 < \ve \le \ve _1 << 1$ and ${\rm dist} (z,W) < \ve$, then
$$
\log \ve ^2 - C_r \le s_{r,\ve} \le 0.
$$
\end{enumerate}
\end{lem}

\begin{proof}
Property 1 is seen as follows:  let $f$ be a test $(n-1,n-1)$-form.  Then
\begin{eqnarray*}
\lim _{\ve \to 0} \int _B (dd^c s_{r,\ve}) \wedge f  &=& \lim _{\ve \to 0} \int _B s_{r,\ve}dd^c f \\
&=&  \int _B s_r dd^c f \\
&=& \int _B ([W] - \Upsilon ^W_r )\wedge f.
\end{eqnarray*}

Property 2 may be established locally, and using group invariance
and uniform flatness, we need only consider the case $z=0$.  But
then by the calculation in the remark following the proof of Lemma
\ref{s-prop} we may assume that $s_r = \log |\zeta ^n|^2$, and
thus 2 follows by integration.
\end{proof}

\subsection{The proof of Theorem \ref{samp-thm}}

\subsubsection*{A positivity lemma}

A key idea behind the proof of the lower sampling inequality is a certain positivity lemma, which we now state and prove.

\begin{lem}\label{pos-lem}
Let $\theta$ be a positive $(n-1,n-1)$-form in $B$ such that for
some weight $\psi$ and each $h \in \sh ^2(B,\psi)$,
$$
\int _B |h|^2 e^{-\psi} \ii \di \dbar \theta < +\infty.
$$
Then
$$
\int _B |h|^2 e^{-\psi} \ii \di \dbar \psi \wedge \theta \ge - \int _B |h|^2 e^{-\psi} \ii \di \dbar \theta.
$$
\end{lem}

\begin{proof}
Letting $S = |h|^2 e^{-\psi}$, one calculates that
$$
\frac{\ii \di \dbar S}{S} = \frac{\ii \di S \wedge \dbar S}{S^2} + \ii \di \dbar \log |h|^2 - \ii \di \dbar \psi,
$$
and thus
$$
\ii \di \dbar S \wedge \theta \ge - S \ii \di \dbar \psi \wedge \theta.
$$

Let $f : \R \to [0,1]$ be a smooth function supported on $(-\infty
, 3/4]$ such that $f |(-\infty, 0] \equiv 1$.  Consider the function
$$
\chi _a (z) = f \left ( \frac{1-|z|^2}{-a} +1 \right ), \qquad a > 0.
$$
Then
\begin{eqnarray*}
\int _B \ii \di \dbar S \wedge \theta &=& \lim _{a \to 0+} \int _B \chi _a(z) \ii \di \dbar S \wedge \theta \\
&=& \lim _{a \to 0+} \int _B S \ii \di \dbar \left ( \chi _a(z)
\wedge \theta \right ) \\
&=& \lim _{a \to 0+} \left ( \int _B \chi _a(z) S \ii \di \dbar \theta + O(a)\right ) \\
&=& \int _B S \ii \di \dbar \theta,
\end{eqnarray*}
where the second equality follows from the Green-Stokes identity (\ref{sg-id}).
\end{proof}

\subsubsection*{Conclusion of the proof of Theorem \ref{samp-thm}}

Let
$$
\psi = \kappa + n \lambda + \alpha s_{r,\ve}.
$$
In view of Lemmas \ref{curb} and \ref{density-conditions}.2, for some $0< <\alpha < 1$, $c > 0$ and $\theta \in \sP _W(B)$ we have
$$
\ii\di \dbar \psi \wedge \theta \le - ce^{n\lambda } \omega _B^n + C[W]_{\ve} \wedge \omega _B ^{n-1},
$$
where $[W]_{\ve}$ denotes the regularization of the current $[W]$ in the manner of Lemma \ref{curb}.  Let $F \in \sh ^2(B,\kappa)$.  Then
\begin{eqnarray*}
\int _B |F|^2 e^{-\kappa} \omega _B ^n &\le & \int _{B}|F|^2 e^{-\psi} \omega _B ^n \\
&\le & C \int _{B} |F|^2 e^{-\psi} [W]_{\ve} \wedge \omega _B ^{n-1} - C \int _{B} |F|^2 e^{-\psi} \ii \di \dbar \psi \wedge \theta\\
&\le &  C \int _{B} |F|^2 e^{-\psi} [W]_{\ve} \wedge \omega _B ^{n-1},
\end{eqnarray*}
where the last inequality follows from Lemma \ref{pos-lem} and the definition of $\sP_W(B)$.  Thus we have
\begin{eqnarray*}
\int _{B} |F|^2 e^{-\kappa } \omega _B ^n & \le &  C \int _{B} |F|^2 e^{-\psi} [W]_{\ve} \wedge \omega _B ^{n-1} \\
&\le & \frac{C}{\ve^2} \int _{N_{\ve}(W)}|F|^2 e^{-\psi} \omega _B ^{n}\\
&\le & \frac{C}{\ve^{2+2\alpha}} \int _{N_{\ve}(W)} |F|^2 e^{-\kappa} \omega _B ^n.\\
\end{eqnarray*}

Our next task is to compare
\[
\int _{N_{\ve}(W)} |F|^2 e^{-\kappa} \omega _B ^n \qquad {\rm with} \qquad \int _{W} |F|^2 e^{-\kappa} \omega _B ^{n-1}.
\]
To do this, we cover $N_{\ve} ^B (W)$ by domains
\[
\Delta _p(\ve)  = \Phi ( (E(p,\ve _p) \cap W) \times \D (0,\ve)), \quad p \in \sw.
\]
Here $\Phi$ is the diffeomorphism defined in the remark in Section \ref{uf-section} following Definition \ref{uni-flat}, $E(p,\ve _p)$ is the Bergman-Green ball of center $p$ and radius $\ve _p$, and $\sw \subset W$ is a discrete set that is uniformly separated with respect to the Bergman-Green distance.  We now employ Lemma \ref{quimbo-trick} once more to obtain a function
\[
\Delta _p (\ve) {\buildrel \Phi \over \cong} ( W \cap E(p,\ve _p))\times \D (0,\ve) \ni (x,t) \mapsto H_p(x,t) \in \C
\]
that is holomorphic in $t$ and satisfies
$$
H_p (x,0)=0 \qquad {\rm and} \qquad \left | 2 \re (H_p(x,t))
+\kappa  (x,0) - \kappa (x,t) \right | \le C
$$
where $C$ is an absolute constant depending only on $\ii \di \dbar \kappa.$

Let $F_p=Fe^{-H_p}$. By Taylor's Theorem, for each $x$ we have
\begin{eqnarray*}
&& |F_p (x,t)|^2 \le  C |F(x,0)|^2 + \ve ^2 \sup _{|t|\le \ve}
\left | \frac{\di F_p}{\di t}\right |^2.
\end{eqnarray*}
We then obtain
\begin{eqnarray*}
&& \int _{\Delta _p(\ve)} |F|^2 e^{-\kappa} \omega _B^n \\
&\le & \int _{\Delta _p(\ve)}|F_p|^2 e^{-\kappa (x,0)} \omega _B ^{n} \\
&\le & C_1\ve ^2 \int _{W\cap E(p,\ve _p)} |F|^2 e^{-\kappa}
\omega _B ^{n-1} + \ve ^2 \int _{\Delta _p(\ve)} \sup _{|t|\le \ve} \left | \frac{\di F _p}{\di t}\right |^2 e^{-\kappa (x,0)}\omega _B^{n} \\
&\le &  C_1 \ve ^2 \int _{W\cap E(p,\ve _p)} |F|^2 e^{-\kappa}
\omega _B ^{n-1} + \ve ^4 \int _{W\cap E(p,\ve _p)} \sup _{|t|\le \ve} \left | \frac{\di F_p}{\di t}\right |^2 e^{-\kappa (x,0)} \omega _B^{n-1}\\
&\le &  C_1 \ve ^2 \int _{W\cap E(p,\ve _p)} |F|^2 e^{-\kappa}
\omega _B^{n-1} + C \ve ^4 \int _{\Delta _p(\ve _o)} |F_p|^2e^{-\kappa (x,0)} \omega _B ^n \\
&\le &  \ve ^2 \int _{W\cap E(p,\ve _p)} |F|^2 e^{-\kappa} \omega ^{n-1} + C \ve ^4 \int _{\Delta_p(\ve _o)} |F|^2 e^{-\kappa }\omega _B ^n,
\end{eqnarray*}
where $\ve < \ve _o /2$ and $\ve_o$ is as in Definition
\ref{uni-flat}. We have used the Cauchy estimates in the
penultimate inequality. Combining all of this, and summing over $p
\in \sw$, we obtain
$$
\int _{B} |F|^2 e^{-\kappa} \omega _B ^n \le \frac{C} {\ve^{2\alpha}} \int _W |F|^2 e^{-\kappa}  \omega _B
^{n-1} + C' \ve ^{2-2\alpha} \int _{B} |F|^2 e^{-\kappa } \omega _B ^n,
$$
which establishes the left inequality in (\ref{samp-ineq}) as soon as we take $\ve$ small enough.  Here we are using the fact that, since $\sw$ is uniformly separated,
\[
\sum _{p \in \sw} \int _{\Delta _p (\ve _o)} |F|^2 e^{-\kappa} \omega ^n _B \le C \int _B |F|^2 e^{-\kappa } \omega ^n _B.
\]
The right inequality was already established in Corollary \ref{upper-sampling-ineq}.  The proof of Theorem \ref{samp-thm} is complete.\qed

\end{document}